\documentclass[12pt]{amsart}

\usepackage{amsmath,amssymb}
\usepackage{amsfonts}
\usepackage{amsthm}
\usepackage{latexsym}
\usepackage{graphicx}
\usepackage{epsfig}

\newtheorem{theorem}{Theorem}[section]
\newtheorem{lemma}[theorem]{Lemma}

\newtheorem{corollary}[theorem]{Corollary}

\newtheorem{prop}{Proposition}[section]
\newtheorem{rema}[prop]{Remark}

\makeatletter \@addtoreset{equation}{section} \makeatother

\def\ddt{\frac{d}{dt}}
\def\ppt{\frac{\partial}{\partial t}}

\begin{document}

\title{Differential Harnack Estimates for
Backward Heat Equations with Potentials under the Ricci Flow}

\author{Xiaodong Cao$^*$
}
\thanks{$^*$Research
partially supported by
%an MSRI postdoctoral fellowship
the Jeffrey Sean Lehman Fund from Cornell University}
% by NSF
%grant no. \# }

\address{Department of Mathematics,
 Cornell University, Ithaca, NY 14853-4201}
\email{cao@math.cornell.edu}
%\curraddr{}
%\thanks{}

%\author{Richard S. Hamilton$^{\flat}$}
%\thanks{$^{\flat}$Research
%partially supported by NSF grant no. \# }

%\address{Department of Mathematics,
 % Columbia University, New York, NY 10027}
%\email{hamilton@math.columbia.edu}
%\curraddr{}
%\thanks{}

%    General info
\renewcommand{\subjclassname}{%
  \textup{2000} Mathematics Subject Classification}
\subjclass[2000]{Primary 53C44}
% Global differential geometry
% 53C21 Methods of Riemannian geometry, including PDE methods;
%   curvature restrictions [See also 58J60]
% 53C44 Geometric evolution equations (mean curvature flow)
% 53C55 Hermitian and K\"ahlerian manifolds [See also 32Cxx]
% Qualitative properties of solutions
% 35B35 Stability, boundedness
%Parabolic equations and systems [See also 35Bxx, 35Dxx, 35R30,
%35R35, 58J35]
%35K55 Nonlinear PDE of parabolic type
%35K90 Abstract parabolic evolution equations
% Partial differential equations on manifolds; differential operators
%58J35 Heat and other parabolic equation methods
% 58J37 Perturbations; asymptotics
% Systems theory; control - Stability
% 93D05 Lyapunov and other classical stabilities
%       (Lagrange, Poisson, $L^p, l^p$, etc.)

%\date{May 18, 2008}

\maketitle

\markboth{Xiaodong Cao} {Differential Harnack Estimates for
Backward Heat Equations with Potentials under the Ricci Flow}

\begin{abstract}  In this paper, we derive a general evolution
formula for possible Harnack quantities. As a consequence, we
prove several differential Harnack inequalities for positive
solutions of backward heat-type equations with potentials
(including the conjugate heat equation) under the Ricci flow. We
shall also derive Perelman's Harnack inequality for the
fundamental solution of the conjugate heat equation under the
Ricci flow.
\end{abstract}

\section{\textbf{Introduction}}
%%%%%%%%%%%%%%%%%%%%%%%%%%%%%%%%%%%%%%%%%%%%%%%%%%%%%%%%%%%%%%%%%%%
In \cite{ly86}, P. Li and S.-T. Yau proved a differential Harnack
inequality by developing a grading estimate for positive solutions
of the heat equation (with fixed metric). More precisely, they
proved that, for any positive solution $f$ of the heat equation
$$\frac{\partial f}{\partial t}=\triangle f$$ on Riemannian
manifolds with nonnegative Ricci curvature, then $$\ppt \ln
f-|\nabla \ln f|^2+\frac{n}{2t}=\triangle \ln f+\frac{n}{2t}\geq
0.$$ The idea was brought to study general geometric evolution
equations by R. Hamilton. The differential Harnack estimate has
since become an important technique in the studies of geometric
evolution equations.

For the Ricci flow, R. Hamilton \cite{harnack} proved  a Harnack
estimate on Riemannian manifolds with weakly positive curvature
operator. The Harnack quantities are also know as curvatures of a
degenerate metric in space-time thanks to the work of B. Chow and
S.-C. Chu \cite{chowchu95}.  In dimension two, R. Hamilton
\cite{Hsurface} proved a Harnack estimate for the scalar curvature
when it is positive, the Harnack estimate in the general case was
proved by B. Chow in \cite{chow91}. B. Chow and R. Hamilton
generalized their results for the heat equation and for the Ricci
flow on surfaces in \cite{chowhamilton97}.

For other geometric flows, R. Hamilton proved  a Harnack estimate
for the mean curvature flow in \cite{mcf95}. B. Chow proved
Harnack estimates for Gaussian curvature flow in \cite{chowgcf91}
and for Yamabe flow in \cite{chowyf92}. H.-D. Cao \cite{caohd92}
proved a Harnack estimate and L. Ni \cite{ni07} proved a matrix
Harnack estimate (of the forward conjugate heat equation) for the
K\"{a}hler-Ricci flow.

For the heat equation, besides the classical result of P. Li and
S.-T. Yau, R. Hamilton proved a matrix Harnack estimate for the
heat equation in \cite{hmatrix93}. C. Guenther \cite{guenther02}
studied the fundamental solution and Harnack inequality of
time-dependent heat equation. In \cite{caoni04}, H.-D. Cao and L.
Ni proved a matrix Harnack estimate for the heat equation on
K\"{a}hler manifolds. In an earlier paper \cite{caohamilton}, R.
Hamilton and the author proved several Harnack estimates for
positive solutions of the heat-type equation with potential when
the metric is evolving by the Ricci flow (a more detailed
discussion about the literature of Harnack estimates can also be
found in that paper). In \cite{perelman1}, G. Perelman proved a
Harnack estimate for the fundamental solution of the conjugate
heat equation under the Ricci flow. Namely, let $(M, g(t))$, $t\in
[0,T]$, be a solution to the Ricci flow on a closed manifold, $f$
be the positive \textit{fundamental} solution to the conjugate
heat equation
$$\ppt f=-\triangle f+Rf,$$ $\tau=T-t$, and $u=-\ln
f-\frac{n}{2}\ln (4\pi \tau)$. Then for $t\in [0,T)$, G. Perelman
proved that
$$2\triangle u-|\nabla u|^2+R+\frac{u}{\tau}-\frac{n}{\tau}\leq
0$$ (see \cite{ni06} or \cite[Chapter 16]{chowetc2} for a detailed
proof).

In the present paper, we will first derive a general evolution
equation for possible Harnack quantities, we will then prove
Harnack estimates for all positive solutions of the backward
heat-type equation with potentials when the metric is evolving
under the Ricci flow.

Suppose $(M,g(t))$, $t\in [0, T]$, is a solution to the Ricci flow
on a closed manifold. Let $f$ be a positive solution of the
backward heat equation with potential $2R$, i.e.,
\begin{align}
  \frac{\partial g_{ij}}{\partial t}=& -2R_{ij}, \label{rf}\\
  \frac{\partial f}{\partial t}=&-\triangle_{g(t)} f+2Rf. \label{heatf2}
\end{align}
Our first main theorem is the following,

\begin{theorem}\label{theorem1.1}
Let $(M, g(t))$, $t\in [0, T]$,
 be a solution
to the Ricci flow on a closed manifold.  Let $f$ be a positive
solution to the backward heat-type equation (\ref{heatf2}),
$u=-\ln f$, $\tau=T-t$ and
$$H=2\triangle u-|\nabla u|^2+2R-\frac{2n}{\tau}.$$ Then for
all time $t\in [0,T)$,
$$H\leq 0.$$
\end{theorem}

If we further assume that our solution to the Ricci flow is of
Type I, i.e.,
$$|Rm|\leq \frac{d_0}{T-t}$$ for some constant $d_0$, here $T$ is
the blow-up time, then we shall prove the following theorem.

\begin{theorem}\label{theorem1.2} Let $(M, g(t))$, $t\in [0, T)$,
 be a type I solution
to the Ricci flow on a closed manifold. Let $f$ be a positive
solution to the backward heat equation (\ref{heatf2}), $u=-\ln f$,
$\tau=T-t$ and
$$H=2 \triangle u-|\nabla u|^2+2R-d\frac{n}{\tau},$$ here
$d=d(d_0,n)\geq 2$ is some constant such that $H(\tau)<0$ for
small $\tau$. Then for all time $t\in [0,T)$,
$$H\leq 0.$$
\end{theorem}

We will consider the conjugate heat equation under the Ricci flow.
In this case, we also assume that our initial metric $g(0)$ has
nonnegative scalar curvature, it is well-known that this property
is preserved by the Ricci flow. We shall prove

\begin{theorem} \label{theorem1.3} Let $(M, g(t))$, $t\in [0, T]$,
 be a solution
to the Ricci flow (\ref{rf}) on a closed manifold, and suppose
that $g(t)$ has nonnegative scalar curvature. Let $f$ be a
positive solution to the conjugate heat equation
\begin{equation}\label{che}
\frac{\partial f}{\partial t}=-\triangle_{g(t)} f+Rf,
\end{equation}
 $u=-\ln f$,
$\tau=T-t$ and
$$H=2 \triangle u-|\nabla u|^2+R-2\frac{n}{\tau}.$$ Then for all
time $t\in [0,T)$,
$$H\leq 0.$$
\end{theorem}

\begin{rema}
S. Kuang and Q. Zhang \cite{kz08} have also proved a similar
estimate as in Theorem \ref{theorem1.3} (see Remark 3.1). Our
proof follows from a direct calculation of more general evolution
equation in Lemma \ref{lemma:H}.
\end{rema}

The rest of this paper is organized as follows. In section 2, we
will first derive a general evolution equation of Harnack quantity
H, for backward heat-type equation with potentials, then we prove
Theorem 1.1 and Theorem 1.2. We will also prove an integral
version of the Harnack inequality (Theorem 2.3). In section 3, we
will prove Theorem 1.3, then we will derive a general evolution
formula for a Harnack quantity similar to Perelman's, as a
consequence, we will prove Perelman's Harnack estimate. In section
4, we will define two entropy functionals and prove that they are
monotone. In section 5, we will prove a gradient estimate for the
backward heat equation (without the potential term). This is also
a consequence of the general evolution formula of our Harnack
quantity.

{\bf Acknowledgement:} The author would like to thank Laurent
Saloff-Coste for general discussion on Harnack inequalities. He
would also like to thank Leonard Gross and Richard Hamilton for
helpful suggestions.

%%%%%%%%%%%%%%%%%%%%%%%%%%%%%%%%%%%%%%%%%%%%%%%%%%%%%%%%%%%%%%%%%%%%%%
\section{\textbf{General Evolution Equation and Proof of Theorem \ref{theorem1.1}}}%%%%%%%%%%%%%%%%%%
%%%%%%%%%%%%%%%%%%%%%%%%%%%%%%%%%%%%%%%%%%%%%%%%%%%%%%%%%%%%%%%%%%%%%%

Let us first consider positive solutions of general evolution
equations
$$\ppt f=-\triangle f
-cRf$$ for all constant $c$ which we will fix later. Let
$f=e^{-u}$, then $\ln f=-u$. We have
$$\ppt \ln f=-\ppt u,$$ and
$$\nabla \ln f=-\nabla u, ~\triangle \ln f=-\triangle u.$$
Hence $u$ satisfies the following equation,
\begin{equation}\label{eqnu}
\ppt u=-\triangle u+|\nabla u|^2+cR.
\end{equation}
Let $\tau=T-t$, then we have
\begin{equation}\label{eqnf2}
\frac{\partial f}{\partial \tau}=\triangle f+cRf,
\end{equation}
and $u$ satisfies
\begin{equation}\label{eqnu2}
\frac{\partial u}{\partial \tau}=\triangle u-|\nabla u|^2-cR.
\end{equation}
We can now define a general Harnack quantity and derive its
evolution equation.

\begin{lemma}\label{lemma:H} Let $(M, g(t))$ be a solution to the Ricci flow, and
$u$ satisfies (\ref{eqnu2}). Let $$H=\alpha \triangle
u-\beta|\nabla u|^2+aR+b\frac{u}{\tau}+d\frac{n}{\tau},$$ where
$\alpha$, $\beta$, $a$, $b$ and $d$ are constants that we will
pick later. Then $H$ satisfies the following evolution equation,
\begin{align*}
\frac{\partial H}{\partial \tau}=&\triangle H-2\nabla H \cdot
\nabla u-(2\alpha-2\beta)|\nabla \nabla
u+\frac{\alpha}{2\alpha-2\beta}R_{ij}-\frac{\lambda}{2\tau}g_{ij}|^2
-\frac{2\alpha-2\beta}{\alpha} \frac{\lambda}{\tau}H
\\&-2(\alpha-2\beta)
R_{ij}u_iu_j+2(a+\beta c) \nabla R \cdot \nabla u
+(2\alpha-2\beta)\frac{n
\lambda^2}{4\tau^2}\\&+(b-\frac{2\alpha-2\beta}{\alpha} \lambda
\beta)\frac{|\nabla u|^2}{\tau}-(\alpha c+2a) \triangle
R+(\frac{\alpha^2}{2\alpha-2\beta}-2a)|Rc|^2\\&+
(\frac{2\alpha-2\beta}{\alpha}
\lambda-1)b\frac{u}{\tau^2}+(\frac{2\alpha-2\beta}{\alpha} \lambda
-1)d\frac{n}{\tau^2}+(\frac{2\alpha-2\beta}{\alpha} a
\lambda-\alpha\lambda -bc)\frac{R}{\tau},
\end{align*}where $\lambda$ is also a constant that we will pick
later.
\end{lemma}

\begin{proof} The proof follows from a direct computation. We
first calculate the first two terms in $H$,
$$\frac{\partial (\triangle u)}{\partial \tau}=\triangle (\triangle u)
-\triangle (|\nabla u|^2)-c\triangle R-2R_{ij}u_{ij},$$ and
\begin{align*}
\frac{\partial (|\nabla u|^2)}{\partial \tau}=&2\nabla u \cdot
\nabla \triangle u-2Rc(\nabla u, \nabla u)-2\nabla u\cdot \nabla
(|\nabla u|^2)-2c\nabla u
\cdot \nabla R\\
=&\triangle (|\nabla u|^2)-2|\nabla \nabla u|^2-2\nabla u \cdot
\nabla (|\nabla u|^2)-2c\nabla u \cdot \nabla R-4R_{ij}u_i u_j,
\end{align*}
here we used $$\triangle (|\nabla u|^2)=2\nabla u \cdot \triangle
\nabla u+2|\nabla \nabla u|^2,$$ and $$\triangle \nabla u=\nabla
\triangle u+Rc(\nabla u,\cdot).$$ Using the evolution equation of
$R$,
$$\frac{\partial R}{\partial \tau}=-\triangle R - 2|Rc|^2,$$
and (\ref{eqnu2}), we have
\begin{align*} \frac{\partial}{\partial \tau} H=&\triangle H
-\alpha \triangle (|\nabla u|^2)-\alpha c \triangle R-2\alpha
R_{ij}u_{ij}+2\beta|\nabla \nabla u|^2+2\beta \nabla u \cdot
\nabla (|\nabla u|^2)\\&+2\beta c\nabla u \cdot \nabla R +4\beta
R_{ij}u_i u_j-2a \triangle R- 2a|Rc|^2-b\frac{|\nabla
u|^2}{\tau}-b\frac{cR}{\tau}-b\frac{u}{\tau^2}-d\frac{n}{\tau^2}\\
=&\triangle H-2\nabla H \cdot \nabla u-2(\alpha-2\beta)
R_{ij}u_iu_j-(2\alpha-2\beta)|\nabla \nabla u|^2+2(a+\beta c)
\nabla R \cdot \nabla u \\&+b\frac{|\nabla u|^2}{\tau}-(\alpha
c+2a) \triangle R- 2\alpha
R_{ij}u_{ij}-2a|Rc|^2-b\frac{u}{\tau^2}-d\frac{n}{\tau^2}
-b\frac{cR}{\tau}\\
=&\triangle H-2\nabla H \cdot \nabla u-2(\alpha-2\beta)
R_{ij}u_iu_j-(2\alpha-2\beta)|\nabla \nabla
u+\frac{\alpha}{2\alpha-2\beta}R_{ij}-\frac{\lambda}{2\tau}g_{ij}|^2
\\&+2(a+\beta c) \nabla R \cdot \nabla u
+(2\alpha-2\beta)\frac{n \lambda^2}{4\tau^2}-(2\alpha-2\beta)
\frac{\lambda}{\tau}\triangle u-\frac{\alpha
\lambda}{\tau}R\\&+b\frac{|\nabla u|^2}{\tau}-(\alpha c+2a)
\triangle
R+(\frac{\alpha^2}{2\alpha-2\beta}-2a)|Rc|^2-b\frac{u}{\tau^2}-d\frac{n}{\tau^2}
-b\frac{cR}{\tau}\\
=&\triangle H-2\nabla H \cdot \nabla u-2(\alpha-2\beta)
R_{ij}u_iu_j-(2\alpha-2\beta)|\nabla \nabla
u+\frac{\alpha}{2\alpha-2\beta}R_{ij}-\frac{\lambda}{2\tau}g_{ij}|^2
\\&-\frac{2\alpha-2\beta}{\alpha}
\frac{\lambda}{\tau}H+2(a+\beta c) \nabla R \cdot \nabla u
+(2\alpha-2\beta)\frac{n
\lambda^2}{4\tau^2}\\&+(b-\frac{2\alpha-2\beta}{\alpha} \lambda
\beta)\frac{|\nabla u|^2}{\tau}-(\alpha c+2a) \triangle
R+(\frac{\alpha^2}{2\alpha-2\beta}-2a)|Rc|^2\\&+
(\frac{2\alpha-2\beta}{\alpha}
\lambda-1)b\frac{u}{\tau^2}+(\frac{2\alpha-2\beta}{\alpha} \lambda
-1)d\frac{n}{\tau^2}+(\frac{2\alpha-2\beta}{\alpha} a
\lambda-\alpha\lambda -bc)\frac{R}{\tau}.
\end{align*}
\end{proof}

In the above lemma, let us take $\alpha=2$, $\beta=1$, $a=2$,
$c=-2$, $\lambda=2$, $b=0$, $d=-2$, as a consequence of the above
lemma, we have
\begin{corollary} Let $(M, g(t))$ be a solution to the
Ricci flow, $f$ be a positive solution of the following backward
heat equation
$$\ppt f=-\triangle f
+2Rf,$$ let $u=-\ln f$, $\tau=T-t$ and $$H=2 \triangle u-|\nabla
u|^2+2R-2\frac{n}{\tau}.$$ Then we have
\begin{align}\label{hgrad1}
\frac{\partial}{\partial \tau} H=&\triangle H -2\nabla H \cdot
\nabla u -\frac{2}{\tau}H-\frac{2}{\tau}|\nabla u|^2-2|Rc|^2-2
|\nabla_i \nabla_j u+R_{ij}-\frac{1}{\tau}g_{ij}|^2.
\end{align}
\end{corollary}

Now we can finish the proof of Theorems \ref{theorem1.1} and 1.2.

\begin{proof} (Proof of Theorems \ref{theorem1.1} and 1.2)
To prove Theorem 1.1, it is easy to see that for $\tau$ small
enough then $H(\tau)<0$.  It follows from (\ref{hgrad1}) and
maximum principle that
$$H\leq 0$$ for all time $\tau$.

From the above proof, we can easily see that if the solution to
the Ricci flow is of type I and $d\geq 2$ is large enough, such
that $H(\tau)<0$ for $\tau$ small, then Theorem 1.2 is true for
all time $t<T$.
\end{proof}

We now can integrate the inequality $$2 \triangle u-|\nabla
u|^2+2R-2\frac{n}{\tau}\leq 0$$ along a space-time path and get a
classical Harnack inequality.

\begin{theorem}\label{thm2:3} Let $(M, g(t))$, $t\in [0, T]$,
 be a solution
to the Ricci flow on a closed manifold. Let $f$ be a positive
solution to the  equation $$\ppt f=-\triangle f +2Rf.$$ Assume
that $(x_1, t_1)$ and $(x_2, t_2)$, $0\leq t_1<t_2<T$, are two
points in $M\times [0, T)$. Let
$$\Gamma=\inf_{\gamma} \int_{t_1}^{t_2} (\frac12 |\dot{\gamma}|^2+R) dt,$$
where $\gamma$ is any space-time path joining $(x_1, t_1)$ and
$(x_2, t_2)$. Then we have
$$f(x_2,t_2)\leq f(x_1, t_1) (\frac{T-t_1}{T-t_2})^n \exp^{\Gamma}.$$
\end{theorem}

\begin{proof}
Since $H\leq 0$, $\tau=T-t$ and $u=-\ln f$ satisfies
$$\frac{\partial u}{\partial \tau}=\triangle u-|\nabla u|^2+2R,$$
we have
$$2\frac{\partial u}{\partial \tau}+|\nabla
u|^2-2R-\frac{2n}{\tau}\leq 0.$$If we pick a space-time path
$\gamma(x,t)$ joining $(x_1,\tau_1)$ and $(x_2,\tau_2)$ with
$\tau_1>\tau_2>0$. Along $\gamma$, we have
\begin{align*}
\frac{d u}{d \tau} &=\frac{\partial u}{\partial
\tau}+\nabla u \cdot \dot{\gamma}\\
&\leq -\frac{1}{2}|\nabla u|^2+R+\frac{n}{\tau}+\nabla u
\cdot \dot{\gamma}\\
&\leq \frac{1}{2}(|\dot{\gamma}|^2+2R)+\frac{n}{\tau}.
\end{align*}
Hence $$u(x_1, \tau_1)-u(x_2,\tau_2)\leq \frac12 \inf_{\gamma}
\int_{\tau_2}^{\tau_1} (|\dot{\gamma}|^2+2R) d\tau +n\ln
(\frac{\tau_1}{\tau_2}).$$ Or we can write this as
$$u(x_1, t_1)-u(x_2,t_2)\leq \frac12 \inf_{\gamma}
\int_{t_1}^{t_2} (|\dot{\gamma}|^2+2R) dt +n\ln
(\frac{T-t_1}{T-t_2}).$$
 If we denote $\Gamma=\inf_{\gamma}
\int_{t_1}^{t_2} (\frac12 |\dot{\gamma}|^2+R) dt$, then we have
$$f(x_2,t_2)\leq f(x_1, t_1) (\frac{T-t_1}{T-t_2})^n
\exp^{\Gamma},$$ this finishes the proof.
\end{proof}

%%%%%%%%%%%%%%%%%%%%%%%%%%%%%%%%%%%%%%%%%%%%%%%%%%%%%%%%%%%%%%%%%%%%%%
\section{\textbf{On the Conjugate Heat Equation}}%%%%%%%%%%%%%%%%%%
%%%%%%%%%%%%%%%%%%%%%%%%%%%%%%%%%%%%%%%%%%%%%%%%%%%%%%%%%%%%%%%%%%%%%%
In this section, we consider positive solutions of general
evolution equations
$$\ppt f=-\triangle f
-cRf$$ on $[0,T]$, for the special case $c=-1$, which is the case
of conjugate heat equation. In Lemma \ref{lemma:H}, let us take
$\alpha=2$, $\beta=1$, $a=1$, $c=-1$, $\lambda=2$, $b=0$, $d=-2$,
and we arrive at
\begin{corollary} Let $(M, g(t))$, $t\in [0,T]$, be a solution to the
Ricci flow, $f$ be a positive solution of the conjugate heat
equation
\begin{equation}\label{che}
\ppt f=-\triangle f +Rf,
\end{equation}
 let $u=-\ln f$, $\tau=T-t$ and
$$H=2 \triangle u-|\nabla u|^2+R-2\frac{n}{\tau}.$$ Then we have
\begin{align}\label{harnackche}
\frac{\partial}{\partial \tau} H=&\triangle H-2\nabla H \cdot
\nabla u-2|u_{ij}+R_{ij}-\frac{1}{\tau}g_{ij}|^2 -\frac{2}{\tau}H
-\frac{2}{\tau}|\nabla u|^2 -2\frac{R}{\tau}.
\end{align}
\end{corollary}

Now we can finish the proof of Theorem \ref{theorem1.3}.

\begin{proof} (Proof of Theorem \ref{theorem1.3}) It is easy to
see that for $\tau$ small enough then $H(\tau)<0$.  It follows
from (\ref{harnackche}) and maximum principle that
$$H\leq 0$$ for all time $\tau$.
\end{proof}

As in section 2, we can see that if the solution to the Ricci flow
is of type I, i.e.,
$$|Rm|\leq \frac{d_0}{T-t},$$ here $T$ is the blow-up time,
then the Harnack estimate is true for all time $t<T$.

\begin{theorem}\label{thm3.2} Let $(M, g(t))$, $t\in [0, T)$,
 be a type I solution
to the Ricci flow on a closed manifold with nonnegative scalar
curvature. Let $f$ be a positive solution to the conjugate heat
equation (\ref{che}), $u=-\ln f$, $\tau=T-t$ and
$$H=2 \triangle u-|\nabla u|^2+R-d\frac{n}{\tau},$$ here
$d=d(d_0,n)\geq 2$ is some constant such that $H(\tau)<0$ for
small $\tau$. Then for all time $t\in [0,T)$,
$$H\leq 0.$$
\end{theorem}

We can also derive a classical Harnack inequality by integrating
along a space-time path.

\begin{theorem}Let $(M, g(t))$, $t\in [0, T]$,
 be a solution
to the Ricci flow on a closed manifold with nonnegative scalar
curvature. Let $f$ be a positive solution to the conjugate heat
equation $$\ppt f=-\triangle f +Rf.$$ Assume that $(x_1, t_1)$ and
$(x_2, t_2)$, $0\leq t_1<t_2<T$, are two points in $M\times [0,
T)$. Let
$$\Gamma=\inf_{\gamma} \int_{t_1}^{t_2} (|\dot{\gamma}|^2+R) dt,$$
where $\gamma$ is any space-time path joining $(x_1, t_1)$ and
$(x_2, t_2)$. Then we have
$$f(x_2,t_2)\leq f(x_1, t_1) (\frac{T-t_1}{T-t_2})^n \exp^{\Gamma/2}.$$
\end{theorem}

In the rest of this section, we will derive Perelman's Harnack
inequality for the conjugate heat equation (\ref{che}). First let
$f=(4\pi \tau)^{-n/2}e^{-v}$, then $\ln f=-\frac{n}{2} \ln (4\pi
\tau)-v$. We have
$$\ppt \ln f=-\ppt v+\frac{n}{2\tau},$$ and
$$\nabla \ln f=-\nabla v, ~\triangle \ln f=-\triangle v.$$
Hence $v$ satisfies the following equation,
\begin{equation}\label{eqnv}
\ppt v=-\triangle u+|\nabla u|^2+cR+\frac{n}{2\tau}.
\end{equation}
Or if we write the equation in backward time $\tau$, we have
\begin{equation}\label{eqnf2}
\frac{\partial f}{\partial \tau}=\triangle f+cRf,
\end{equation}
and $v$ satisfies
\begin{equation}\label{eqnv2}
\frac{\partial v}{\partial \tau}=\triangle v-|\nabla
v|^2-cR-\frac{n}{2\tau}.
\end{equation}
We can now define a general Harnack quantity and derive its
evolution equation.

\begin{lemma} \label{lemmaph}Let $(M, g(t))$, $t\in [0,T]$, be a solution to the
Ricci flow, and
$v$ satisfies (\ref{eqnv2}). Let $$P=\alpha \triangle
v-\beta|\nabla v|^2+aR+b\frac{v}{\tau}+d\frac{n}{\tau},$$ where
$\alpha$, $\beta$, $a$, $b$ and $d$ are constants that we will
pick later. Then $P$ satisfies the following evolution equation,
\begin{align*}
\frac{\partial P}{\partial \tau}=&\triangle P-2\nabla P \cdot
\nabla v-(2\alpha-2\beta)|\nabla \nabla
v+\frac{\alpha}{2\alpha-2\beta}R_{ij}-\frac{\lambda}{2\tau}g_{ij}|^2
-\frac{2\alpha-2\beta}{\alpha} \frac{\lambda}{\tau}P
\\&-2(\alpha-2\beta)
R_{ij}v_iv_j+2(a+\beta c) \nabla R \cdot \nabla v
+(2\alpha-2\beta)\frac{n
\lambda^2}{4\tau^2}\\&+(b-\frac{2\alpha-2\beta}{\alpha} \lambda
\beta)\frac{|\nabla v|^2}{\tau}-(\alpha c+2a) \triangle
R+(\frac{\alpha^2}{2\alpha-2\beta}-2a)|Rc|^2\\&+
(\frac{2\alpha-2\beta}{\alpha}
\lambda-1)b\frac{v}{\tau^2}+(\frac{2\alpha-2\beta}{\alpha} \lambda
-1)d\frac{n}{\tau^2}-\frac{bn}{2\tau^2}+(\frac{2\alpha-2\beta}{\alpha}
a \lambda-\alpha\lambda -bc)\frac{R}{\tau},
\end{align*} here $\lambda$ is also a constant that we will pick
later.
\end{lemma}

\begin{proof} The proof again follows from the same direct computation
as in the proof of Lemma 2.1. Notice that the only extra term
$-\frac{bn}{2\tau^2}$ comes from the evolution of $v$.
\end{proof}

In the above lemma, let us take $\alpha=2$, $\beta=1$, $a=1$,
$c=-1$, $\lambda=1$, $b=1$, $d=-1$, as a consequence of Lemma
\ref{lemmaph}, we have Perelman's Harnack inequality.
\begin{theorem} (Perelman) Let $(M, g(t))$, $t\in [0,T]$, be a
solution to the Ricci flow, $f$ be the  positive fundamental
solution of the conjugate heat equation
$$\ppt f=-\triangle f
+Rf,$$ let $v=-\ln f-\frac{n}{2}\ln (4\pi \tau)$, $\tau=T-t$ and
$$P=2 \triangle v-|\nabla v|^2+R+\frac{v}{\tau}-\frac{n}{\tau}.$$ Then we have
\begin{align}\label{harnack}
\frac{\partial}{\partial \tau} P=&\triangle P-2\nabla P \cdot
\nabla v-2|v_{ij}+R_{ij}-\frac{1}{2\tau}g_{ij}|^2
-\frac{1}{\tau}P.
\end{align}
Moreover, $$P\leq 0$$ on $[0,T)$.
\end{theorem}

In the same spirit of searching Perelman's Harnack inequality, we
shall also take $\alpha=2$, $\beta=1$, $a=1$, $c=-1$, $\lambda=2$,
$b=0$, $d=-2$ in the above Lemma \ref{lemmaph}. We have a Harnack
inequality for all positive solutions of the conjugate heat
equation (\ref{che}).
\begin{theorem} \label{thmnph} Let $(M, g(t))$, $t\in [0,T]$, be a solution to the
Ricci flow, suppose that $g(t)$ has nonnegative scalar curvature.
Let $f$ be a positive solution of the conjugate heat equation
$$\ppt f=-\triangle f
+Rf,$$ let $v=-\ln f-\frac{n}{2}\ln (4\pi \tau)$, $\tau=T-t$ and
$$P=2 \triangle v-|\nabla v|^2+R-\frac{2n}{\tau}.$$ Then we have
\begin{align}\label{harnack}
\frac{\partial}{\partial \tau} P=&\triangle P-2\nabla P \cdot
\nabla v-2|v_{ij}+R_{ij}-\frac{1}{\tau}g_{ij}|^2
-\frac{2}{\tau}P-2\frac{|\nabla v|^2}{\tau}-2\frac{R}{\tau}.
\end{align}Moreover, for all time $t\in [0,T)$,
$$P\leq 0.$$
\end{theorem}

\begin{proof}  It is easy to
see that for $\tau$ small enough then $P(\tau)<0$.  It follows
from (\ref{harnack}) and maximum principle that
$$P\leq 0$$ for all time $\tau$, hence for all $t$.
\end{proof}

\begin{rema}
Theorem \ref{thmnph} can be deduced from Theorem \ref{theorem1.3}
directly since $v=u-\frac{n}{2} \ln (4\pi \tau)$ and $b=0$. But as
we can see in the direct calculation, it will also lead to
Perelman's Harnack inequality (by choosing different
coefficients). S. Kuang and Q. Zhang proved this estimate in
\cite{kz08}.
\end{rema}

\begin{rema}
We can prove similar estimate for $P$ if we have a Type I solution
to the Ricci flow. We can also prove a classical Harnack
inequality by integrating along a space-time path.
\end{rema}

\section{\textbf{Entropy Formulas and Monotonicities}}

In this section, we will define two entropies which are similar to
Perelman's entropy functionals as in \cite{perelman1}, and we will
show that both of them are monotone under the Ricci flow. Let
$(M,g(t))$ be a solution to the Ricci flow on a close manifold, we
shall also assume that $g(t)$ has nonnegative scalar curvature,
 we first prove

\begin{theorem}\label{thm4.1} Assume that $(M, g(t))$, $t\in [0,T]$,
 is a solution to the Ricci
flow on a Riemannian manifold with nonnegative scalar curvature.
Let $f$ be a positive solution of $$\ppt f=-\triangle f +2Rf,$$
$u=-\ln f$ and $\tau=T-t$. Define
$$H=2 \triangle u-|\nabla u|^2+2R-2\frac{n}{\tau}$$ and
$$F=\int_M \tau^2He^{-u} d\mu,$$ then $\forall t\in [0,T)$, we have
$F\leq 0$ and
$$\ddt F\geq 0.$$
\end{theorem}

\begin{proof} The fact that $F\leq 0$ follows directly from $H\leq
0$. We calculate its time derivative, using (\ref{hgrad1}) in
Lemma \ref{lemma:H} and $\ppt d\mu=-Rd\mu$, we have
\begin{align*}\ddt F=-\frac{d F}{d \tau}=&-\int_M (2\tau He^{-u}+\tau^2e^{-u}\ppt
H+\tau^2H\frac{\partial}{\partial \tau} e^{-u}+R\tau^2He^{-u}) d\mu\\
=&-\int_M [\triangle (\tau^2He^{-u})-2\tau^2e^{-u}
|u_{ij}+R_{ij}-\frac{1}{t}g_{ij}|^2-2\tau e^{-u}|\nabla
u|^2\\
&-\tau^2e^{-u}(R +2|Rc|^2)] d\mu \geq 0.
\end{align*}
\end{proof}

We now define an entropy associate with the conjugate heat
equation.

\begin{theorem}\label{thm4.2} Assume that $(M, g(t))$, $t\in [0,T]$,
is a solution to the Ricci flow on a closed Riemannian manifold
with nonnegative scalar curvature. Let $f$ be a positive solution
of $$\ppt f=-\triangle f +Rf,$$ $u=-\ln f$ and $\tau=T-t$. Let
$$H=2 \triangle u-|\nabla u|^2+R-2\frac{n}{\tau}$$ and
$$W=\int_M \tau^2He^{-u} d\mu,$$ then $\forall t\in [0,T)$, we have
$W\leq 0$ and
$$\ddt W\geq 0.$$
\end{theorem}

\begin{proof} The fact that $W\leq 0$ follows directly from $H\leq
0$. To calculate its time derivative, using (\ref{harnackche}) and
$\ppt d\mu=-Rd\mu$, we have
\begin{align*}\ddt W=-\frac{d W}{d \tau}=&-\int_M (2\tau He^{-u}+
\tau^2e^{-u}\frac{\partial}{\partial \tau} H
+\tau^2H\frac{\partial}{\partial \tau} e^{-u}+R\tau^2He^{-u}) d\mu\\
=&-\int_M [\triangle (\tau^2 He^{-u})-2\tau^2e^{-u}
|u_{ij}+R_{ij}-\frac{1}{t}g_{ij}|^2-2\tau e^{-u}|\nabla
u|^2\\
&-2\tau R e^{-u}] d\mu \geq 0.
\end{align*}
\end{proof}

\begin{rema}
As noted in \cite{perelman1}, if we consider the system
$$\left\{\begin{array}{l} \ppt g_{ij}=-2(R_{ij}+\nabla_i \nabla_j
u),\\ \ppt u=-\triangle u-R,
\end{array}\right.$$ then the measure $dm=e^{-u}d\mu$ is fixed.
This system differs from the original system of the conjugate heat
equation under the Ricci flow
$$\left\{\begin{array}{l} \ppt g_{ij}=-2R_{ij},\\ \ppt u=-\triangle
u+|\nabla u|^2-R
\end{array}\right.$$ by a diffeomorphism.
\end{rema}

\section{\textbf{Gradient Estimate for the Backward Heat Equation}}

In this section, we consider a gradient estimate for positive
solutions $f$ to the backward heat equation
\begin{equation}\label{heatf}
\ppt f=-\triangle f.
\end{equation}
Since the equation is linear, without loss of generality, we may
assume that $0<f< 1$. Let $f=e^{-u}$, then $u$ satisfies
\begin{equation}\label{heatu}
\ppt u=-\triangle u+|\nabla u|^2,
\end{equation} and $u> 0$.

In the proof of Lemma \ref{lemma:H}, let take $\alpha=0$,
$\beta=-1$, $a=c=0$, $b=-1$ and $d=0$, then
$$H=|\nabla u|^2-\frac{u}{\tau},$$ and we have
\begin{align}\label{hgrad}
\frac{\partial}{\partial \tau} H=&\triangle H -2\nabla H \cdot
\nabla u-4Rc(\nabla u, \nabla u) -\frac{1}{\tau}H-2 |\nabla \nabla
u|^2.
\end{align}
This leads to the following theorem.
\begin{theorem}
Let $(M, g(t))$, $t\in [0, T]$,
 be a solution
to the Ricci flow on a closed manifold with nonnegative curvature
operator. Let $f(< 1)$ be a positive solution to the backward heat
equation (\ref{heatf}), $u=-\ln f$, $\tau=T-t$ and
$$H=|\nabla u|^2-\frac{u}{\tau}.$$ Then for
all time $t\in [0,T)$,
$$H\leq 0,$$ i.e.,
$$ |\nabla f|^2 \leq \frac{f^2 \ln (1/f)}{T-t}.$$
\end{theorem}

\begin{proof}
Notice that as $\tau$ small enough, $H<0$, now the proof follows
from (\ref{hgrad}) and the maximum principle.
\end{proof}

\bibliographystyle{halpha}
\bibliography{bio}
\end{document}